\documentclass{amsart}
\usepackage{amsmath}
\usepackage{array}
\usepackage{amsfonts}
\usepackage{latexsym}
\usepackage{epsfig}
\newtheorem{theorem}{Theorem}
\newtheorem{prop}[theorem]{Proposition}

\newtheorem{lemma}[theorem]{Lemma}

\newcommand\beq{\begin{equation}}
\newcommand\eeq{\end{equation}}
\newcommand\bce{\begin{center}}
\newcommand\ece{\end{center}}
\newcommand\bea{\begin{eqnarray}}
\newcommand\eea{\end{eqnarray}}
\newcommand\ben{\begin{enumerate}}
\newcommand\een{\end{enumerate}}

\newcommand\wt{\widetilde}

\newcommand\nn{\nonumber}

\newcommand\ms{\medskip}
\newcommand\wh{\widehat}
\newcommand\brr{\begin{array}}
\newcommand\err{\end{array}}
\newcommand\bt{\begin{tabular}}
\newcommand\et{\end{tabular}}
\newcommand\bit{\begin{itemize}}
\newcommand\eit{\end{itemize}}

\newcommand\D{{\mathcal D}}
\newcommand\DD{{\mathcal D}^{(2)}}
\newcommand\T{{\mathcal T}}
\newcommand\TT{{\mathcal T}^{(2)}}
\newcommand\Tk{{\mathcal T}^{(k)}}
\newcommand\La{\Lambda}

\newcommand\bij{\Psi}

\newenvironment{abstrac}{%
         \small
        \begin{center}%
          {\bfseries {Abstract}\vspace{-.5em}}%
                 \end{center}%
        \quotation}

\title{A bijection between $2$-triangulations and pairs of non-crossing Dyck paths}
\author{Sergi Elizalde}
\address{Department
of Mathematics, Dartmouth College, Hanover NH
02139.}\email{sergi.elizalde@dartmouth.edu}

\begin{document}
 \maketitle
 \thispagestyle{empty}

\begin{abstrac}
\medskip

A $k$-triangulation of a convex polygon is a maximal set of
diagonals so that no $k+1$ of them mutually cross in their
interiors. We present a bijection between $2$-triangulations of a
convex $n$-gon and pairs of non-crossing Dyck paths of length
$2(n-4)$. This solves the problem of finding a bijective proof of a
result of Jonsson for the case $k=2$. We obtain the bijection by
constructing isomorphic generating trees for the sets of
$2$-triangulations and pairs of non-crossing Dyck paths.
\end{abstrac}

\section{Introduction}

A triangulation of a convex $n$-gon can be defined as a maximal set
of diagonals so that no two of them intersect in their interiors. It
is well known that the number of triangulations of a convex $n$-gon
is the Catalan number $C_{n-2}=\frac{1}{n-1}\binom{2(n-2)}{n-2}$,
and that all such triangulations have $n-3$ diagonals (not counting
the $n$ sides of the polygon as diagonals).

We say that two diagonals \emph{cross} if they intersect in their
interiors. Define a $m$-crossing to be a set of $m$ diagonals where
any two of them mutually cross. A natural way to generalize a
triangulation is to allow diagonals to cross, but to forbid
$m$-crossings for some fixed $m$. For any positive integer $k$,
define a \emph{$k$-triangulation} to be a maximal set of diagonals
not containing any $(k+1)$-crossing. For example, a
$1$-triangulation is just a triangulation in the standard sense.
Generalized triangulations appear in \cite{CP,DKM,DKM2,Jon,Nak}. It
was shown in \cite{DKM,Nak} that all $k$-triangulations of a convex
$n$-gon have the same number of diagonals. Counting also the $n$
sides of the polygon, the total number of diagonals and sides in a
$k$-triangulation is always $k(2n-2k-1)$.

Jacob Jonsson \cite{Jon} enumerated $k$-triangulations of a convex
$n$-gon, proving the following remarkable result.

\begin{theorem}\label{th:Jon}
The number of $k$-triangulations of a convex $n$-gon is equal to the
determinant
\beq\label{eq:det} \det(C_{n-i-j})_{i,j=1}^k=\left|\begin{array}{ccccc}
C_{n-2} & C_{n-3} & \ldots & C_{n-k} & C_{n-k-1} \\
C_{n-3} & C_{n-4} & \ldots & C_{n-k-1} & C_{n-k-2} \\ \vdots &
\vdots & \ddots & \vdots & \vdots \\ C_{n-k-1} & C_{n-k-2} & \ldots
& C_{n-2k+1} & C_{n-2k} \end{array} \right|,\eeq where
$C_m=\frac{1}{m+1}\binom{2m}{m}$ is the $m$-th Catalan number.
\end{theorem}

On the other hand, it can be shown \cite{DCV} using the lattice path
determinant formula of Lindstr\"om \cite{Lin}, Gessel and Viennot
\cite{GV} that this determinant counts certain fans of non-crossing
lattice paths. Indeed, recall that \emph{Dyck path} can be defined
as a lattice path with north steps $N=(0,1)$ and east steps
$E=(1,0)$ from the origin $(0,0)$ to a point $(m,m)$, with the
property that it never goes below the diagonal $y=x$. We say that
$m$ is the {\em size} or {\em semilength} of the path. The number of
$k$-tuples $(P_1,P_2,\ldots,P_k)$ of Dyck paths from $(0,0)$ to
$(n-2k,n-2k)$ such that each $P_i$ never goes below $P_{i+1}$ is
given by the same determinant (\ref{eq:det}).

In the case $k=1$, this determinant is just $C_{n-2}$, which counts
Dyck paths from $(0,0)$ to $(n-2,n-2)$. There are several simple
bijections between triangulations of a convex $n$-gon and such paths
(see for example \cite[Problem 6.19]{EC2}). However, for $k\ge2$,
the problem becomes more complicated. One of the main open questions
left in \cite{Jon}, stated also in \cite[Problem 1]{Kra}, is to find
a bijection between $k$-triangulations and $k$-tuples of
non-crossing Dyck paths, for general $k$. In this paper we solve
this problem for $k=2$, that is, we find a bijection between
$2$-triangulations of a convex $n$-gon and pairs $(P,Q)$ of Dyck
paths from $(0,0)$ to $(n-4,n-4)$ so that $P$ never goes below $Q$.

In Section~\ref{sec:bij} we present the bijection explicitly. In
Section~\ref{sec:2triang} we describe a generating tree for
$2$-triangulations, and in Section~\ref{sec:Dyck} we give a
generating tree for pairs of non-crossing Dyck paths. In
Section~\ref{sec:why} we show that these two generating trees are
isomorphic, and that our bijection maps each node of one tree to the
corresponding node in the other. In Section~\ref{sec:k} we discuss
possible generalizations of our results to arbitrary $k$.

\subsection{Notation}

From now on, the term $n$-gon will refer to a convex $n$-gon, which
can be assumed to be regular. We label its vertices clockwise with
the integers from $1$ to $n$. For any $n>2k>0$, let $\T^{(k)}_n$
denote the set of $k$-triangulations of an $n$-gon. Let $\D^{(k)}_m$
denote the set of $k$-tuples $(P_1,P_2,\ldots,P_k)$ of Dyck paths
from $(0,0)$ to $(m,m)$ such that $P_i$ never goes below $P_{i+1}$
for $1\le i\le k-1$.

Given $n$ points labeled $1,2,\ldots,n$, a segment connecting $a$
and $b$ (with $a<b$) can be associated to the square $(a,b)$ in an
$n\times n$ board with rows indexed increasingly from top to bottom
and columns from left to right. A collection of segments connecting
some of the points can then be represented as a subset of the
squares of the triangular array $\Omega_n=\{(a,b):1\le a<b\le n\}$,
as it was done in \cite{Jon}. If the points are the vertices of an
$n$-gon labeled clockwise, then the squares $(a,a+1)$, for $1\le
a\le n-1$, and $(1,n)$ correspond to the sides of the polygon. The
remaining squares of $\Omega_n$ correspond to diagonals. The
diagonal connecting two vertices $a$ and $b$ will be denoted
$(a,b)$.

It is easy to check (see for example \cite{Jon}) that $t$ diagonals
$(a_1,b_1),\ldots,(a_t,b_t)$ with $a_1\le a_2\le \cdots\le a_t$ and
$a_i<b_i$ for all $i$ form a $t$-crossing if and only if
$a_1<a_2<\cdots<a_t<b_1<b_2<\cdots<b_t$. The condition that
$a_t<b_1$ can be replaced with the condition that smallest rectangle
containing the $t$ squares $(a_i,b_i)$, $1\le i\le t$, fits inside
$\Omega_n$.

\begin{figure}[hbt]
\epsfig{file=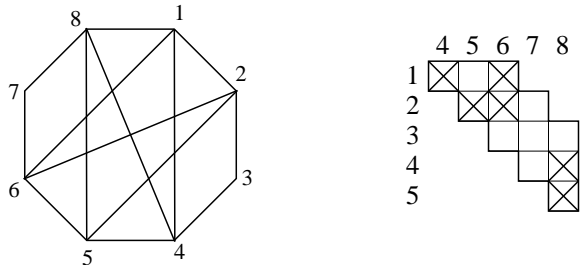,width=3in} \caption{\label{fig:2triang} A
$2$-triangulation of an octagon and its representation as a subset
of $\Lambda_8$.}
\end{figure}

Note that the diagonals joining two vertices that have less than $k$
vertices in between them can never be part of a $k$-crossing. We
will call these {\em trivial diagonals}. They are those of the form
$(a,a+j)$ (or $(a+j-n,a)$ if $a+j>n$), for $2\le j\le k$, $1\le a\le
n$. Any $k$-triangulation of the polygon contains all these
diagonals. For simplicity, we will ignore trivial diagonals.
Deleting from $\Omega_n$ the squares corresponding to trivial
diagonals and to the sides of the polygon, we get the shape
$\Lambda^{(k)}_n=\{(a,b):\ 1\le a<b-k\le n-k,\ a>b-n+k\}$. We will
represent $k$-triangulations as subsets of the squares of
$\Lambda^{(k)}_n$. We will draw a cross in a square to indicate that
the corresponding diagonal belongs to the $k$-triangulation. The
number of crosses is then precisely $k(n-2k-1)$, since that is the
number of diagonals of a $k$-triangulation after the superfluous
ones have been omitted \cite{DKM}. See Figure~\ref{fig:2triang} for
an example of a $2$-triangulation of an octagon, where the trivial
diagonals have been omitted. To simplify notation, $\Lambda^{(2)}_n$
will be denoted $\Lambda_n$.

\section{The bijection} \label{sec:bij}

In this section we give a bijection $\bij$ between
$2$-triangulations of an $n$-gon and pairs $(P,Q)$ of Dyck paths
from $(0,0)$ to $(n-4,n-4)$ so that $P$ never goes below $Q$. We
assume that $n\ge5$.

Let $T\in\TT_n$ be a $2$-triangulation of an $n$-gon. The number of
diagonals, not counting the trivial ones (which are present in any
$2$-triangulation) is $2n-10$. We represent $T$ by placing $2n-10$
crosses in $\Lambda_n$. Index the columns of $\La_n$ from $4$ to
$n$, so that the leftmost column is called ``column $4$", and index
the rows from $1$ to $n-3$. This way, a cross in row $a$ and column
$b$ corresponds to the diagonal $(a,b)$.

\begin{figure}[hbt]
\epsfig{file=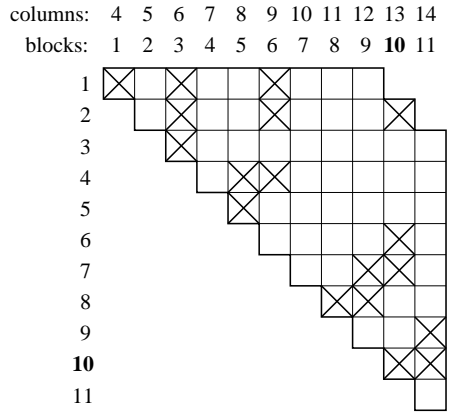,width=2.3in}
\caption{\label{fig:big2triang} A $2$-triangulation of a $14$-gon,
with $r=10$.}
\end{figure}

In the first part of the bijection we will color half of these
crosses blue and the other half red. Along the process, some
adjacent columns of $\La_n$ will be merged. We use the term {\em
block} to refer to a column or to a set of adjacent columns that
have been merged. Blocks are ordered from left to right, so that
``block $j$" refers to the one that has $j-1$ blocks to its left.
At the beginning there are $n-3$ blocks, and block $j$ contains only
column $j+3$, for $j=1,\ldots,n-3$ (see
Figure~\ref{fig:big2triang}). Next we describe an iterative step
that will be repeated $n-5$ times. At each iteration one cross will
be colored blue, another one red, and two blocks will be merged into
one. At the end, all $2n-10$ crosses will be colored, and there will
be only $2$ blocks.

Here is the part that is iterated:
\begin{itemize}
\item Let $r$ be the largest index so that row $r$ has a cross in block $r$.
\item Color blue the leftmost uncolored cross in block $r$ (in case of a tie pick one, for example the
lowest one).
\item Merge blocks $r-2$ and $r-1$ (if $r=2$, we
consider that block $1$ disappears when it is merged with ``block
0").
\item Color red the rightmost uncolored cross in the merged block (in case of a tie pick one, for example the
highest one).
\end{itemize}

\ms

\begin{figure}[hbt]
\bt{ccc}\epsfig{file=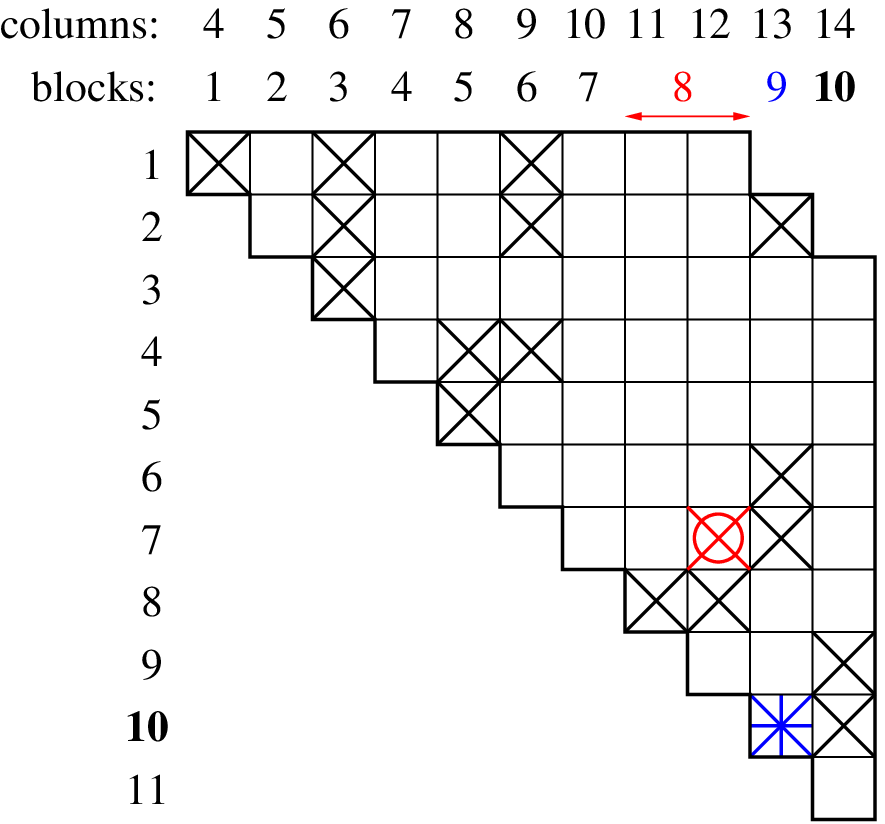,width=1.85in}&\epsfig{file=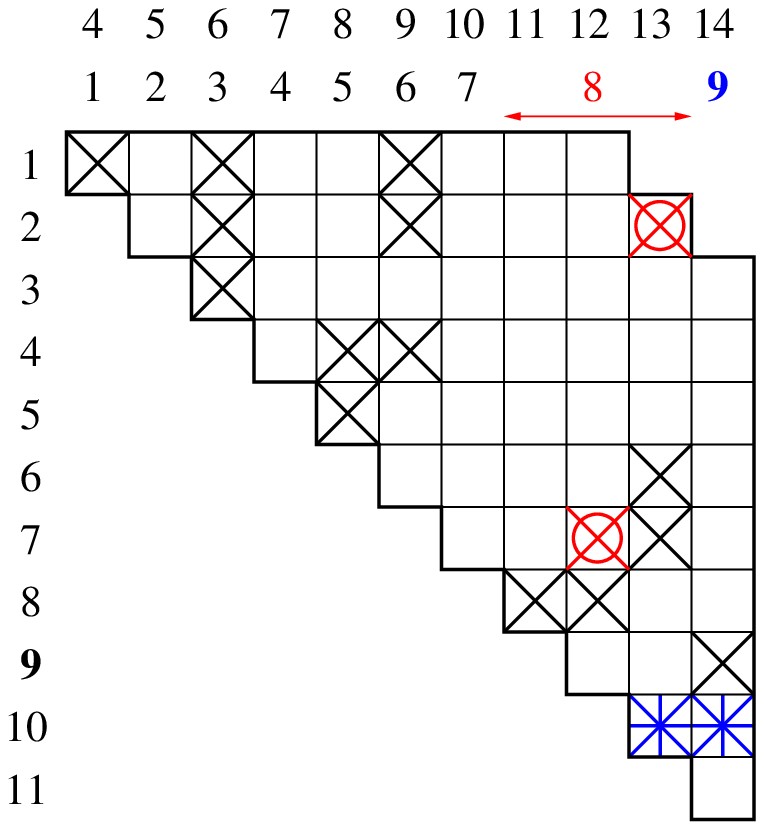,width=1.85in}&\epsfig{file=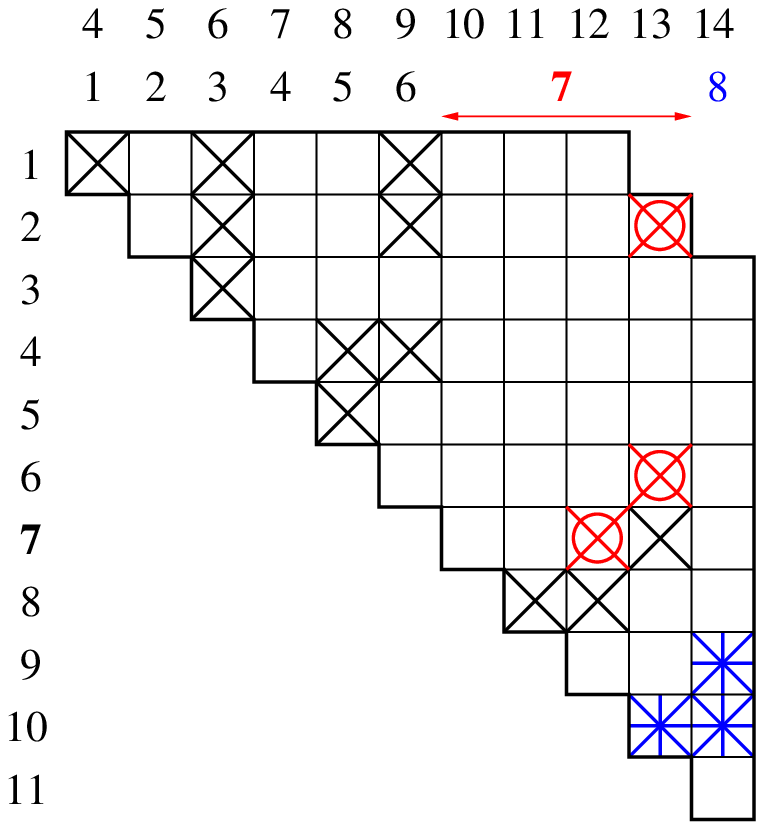,width=1.85in}\\
(a) & (b) & (c) \vspace{5mm}\\
\epsfig{file=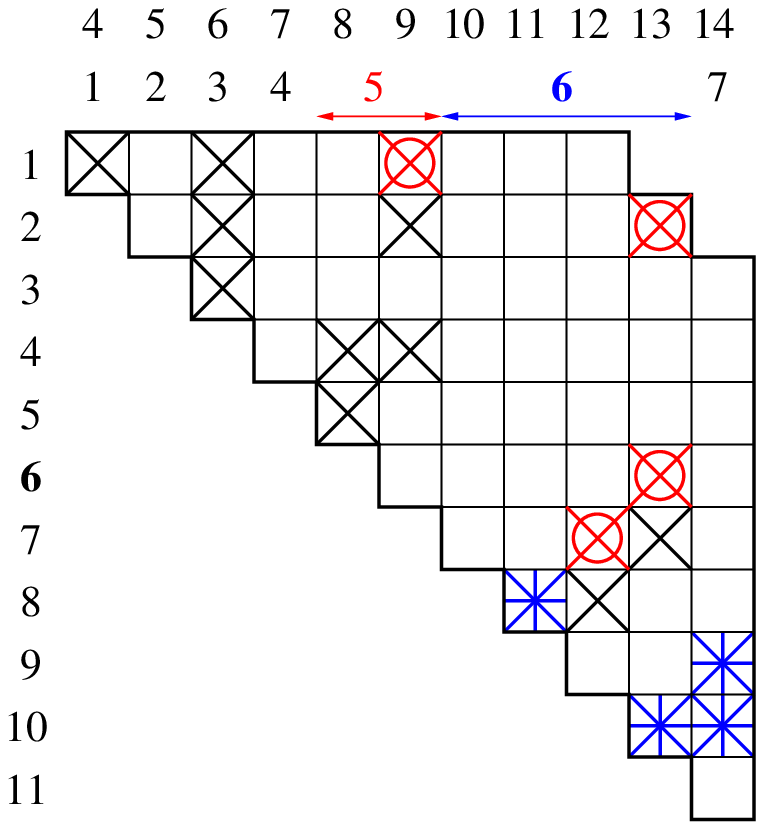,width=1.85in}&\epsfig{file=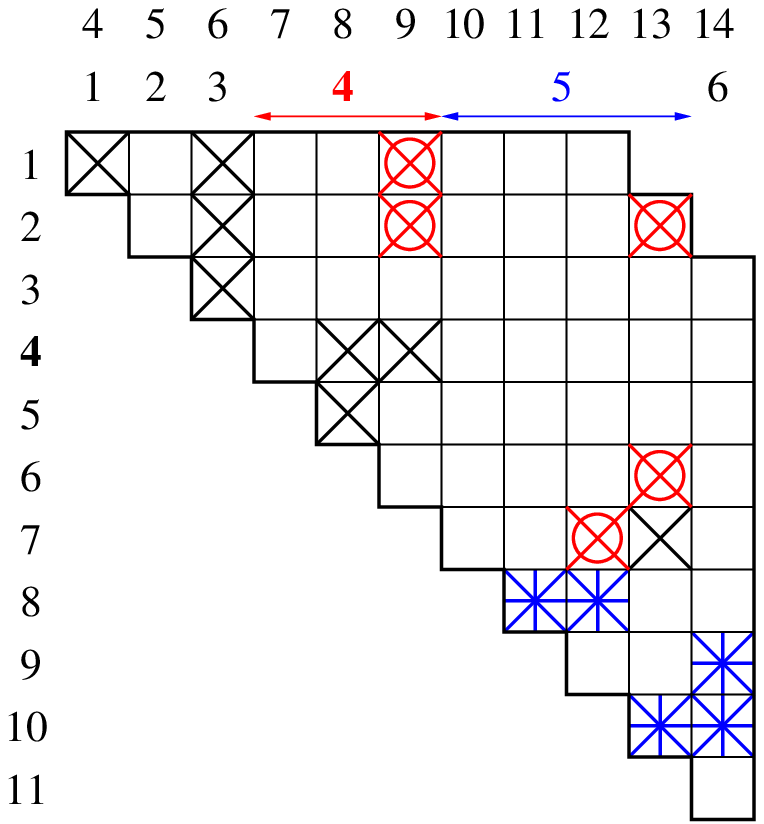,width=1.85in}&\epsfig{file=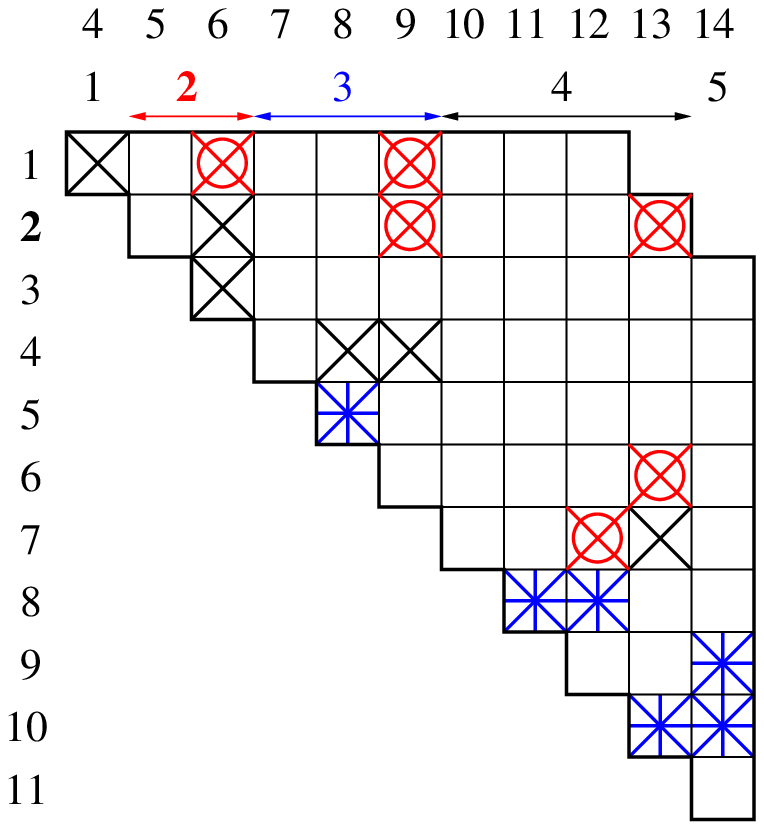,width=1.85in}\\
(d) & (e) & (f) \vspace{5mm}\\
\epsfig{file=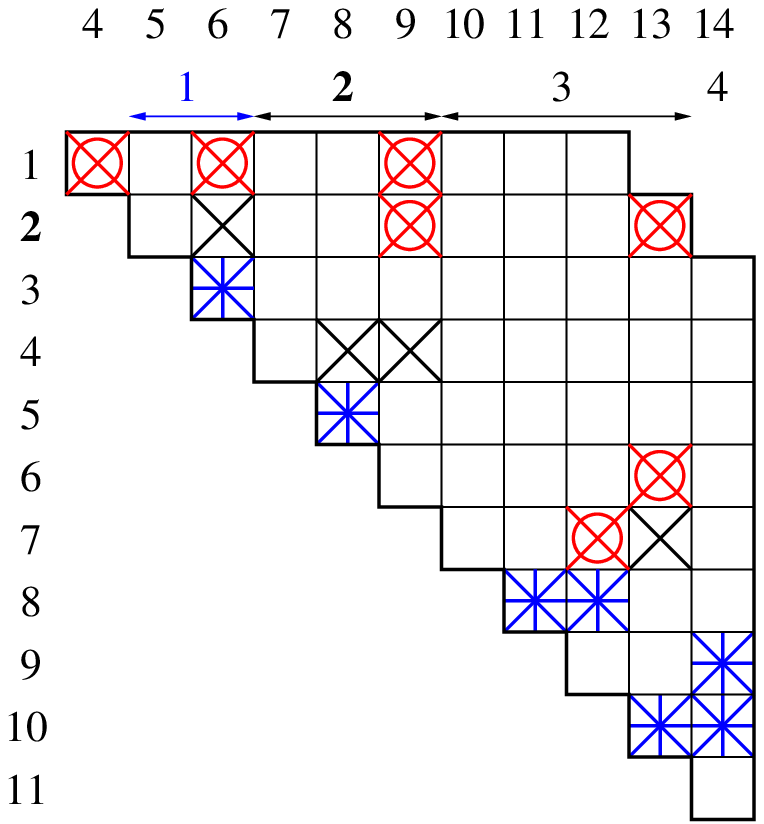,width=1.85in}&\epsfig{file=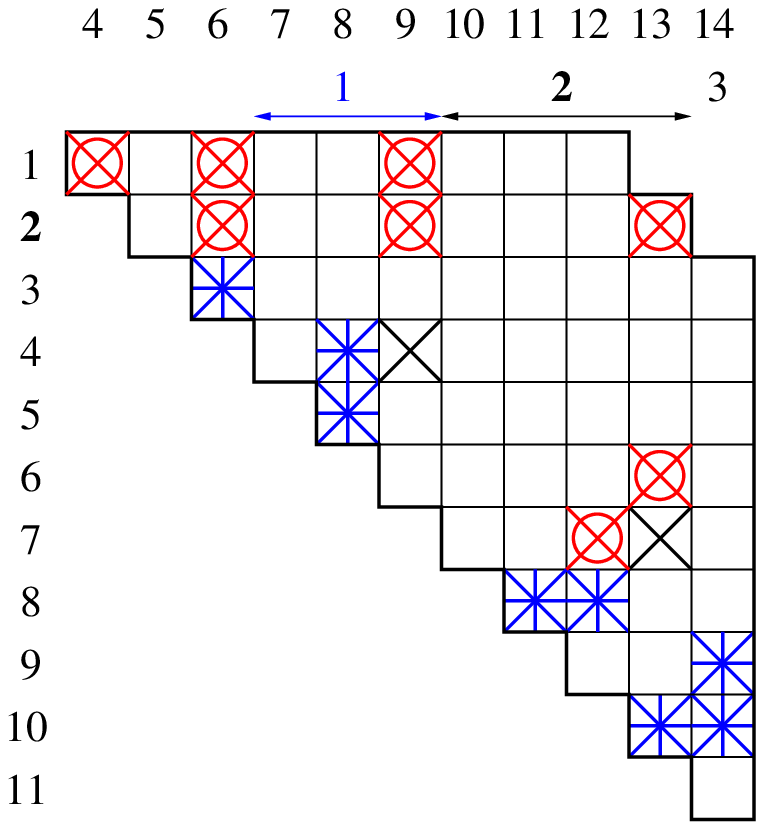,width=1.85in}&\epsfig{file=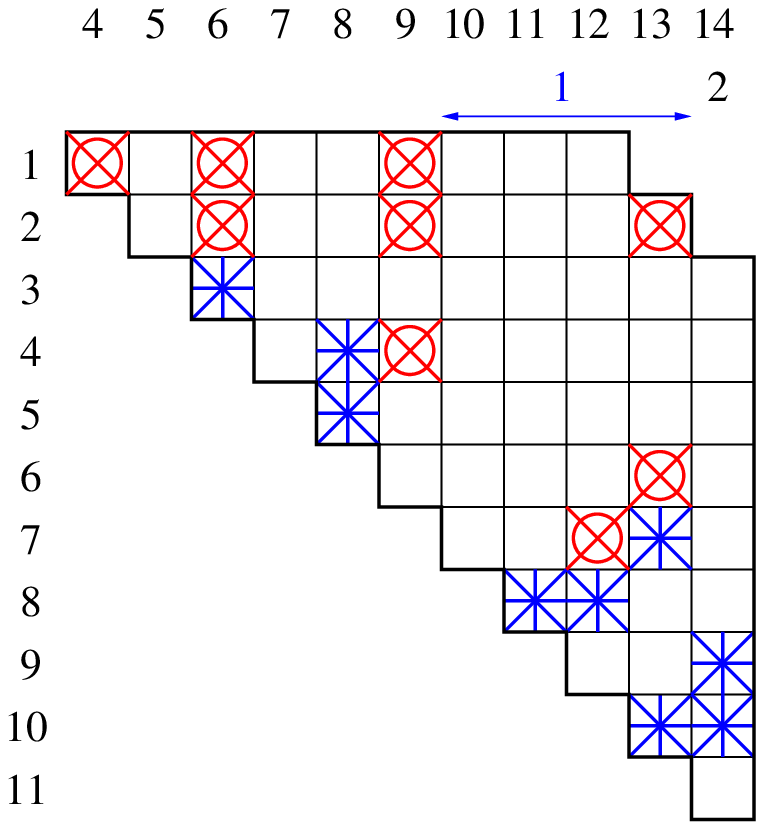,width=1.85in}\\
(g) & (h) & (i) \et\caption{\label{fig:big2triang1} An example of
the coloring algorithm.}
\end{figure}

Let us see how crosses are colored in a particular example. Consider
the $2$-triangulation of a $14$-gon shown in
Figure~\ref{fig:big2triang}. In the following pictures, red crosses
will be drawn with a circle around them, and blue crosses will be
drawn as a star. At the beginning there are 11 blocks, and $r=10$.
In the first iteration, a cross in column 10 is colored blue, a
cross in column 9 is colored red, and columns 8 and 9 are merged
into one block, leaving us with Figure~\ref{fig:big2triang1}(a). In
the second iteration, we have again $r=10$. A cross in block 10 is
colored blue, blocks 8 and 9 are merged, and the leftmost uncolored
cross in the merged block is colored red, as shown in
Figure~\ref{fig:big2triang1}(b). In the third iteration, $r=9$, and
we get Figure~\ref{fig:big2triang1}(c). In the fourth iteration,
$r=7$, so blocks 5 and 6 are merged, giving
Figure~\ref{fig:big2triang1}(d). Next, $r=6$, and blocks 4 and 5 are
merged. In the sixth iteration, $r=4$, and we get
Figure~\ref{fig:big2triang1}(f). In the next step, $r=2$, so block 1
disappears and the cross that it contained is colored red (see
Figure~\ref{fig:big2triang1}(g)). In the last two iterations, $r=2$
again, and we end with Figure~\ref{fig:big2triang1}(i), where all
the crosses have been colored.

\ms

In the second part of the bijection we construct a pair of
non-crossing Dyck paths out of the colored diagram of crosses. For
$j=4,\ldots,n$, let $\alpha_j$ (resp. $\beta_j$) be the number
of blue (resp. red) crosses in column $j$ of $\La_n$. 
Let \bea P&=&NE^{\alpha_5}NE^{\alpha_6}\cdots
NE^{\alpha_{n-1}}NE^{\alpha_{n}}E, \nn \\ \nn
Q&=&NE^{\beta_4}NE^{\beta_5}\cdots
NE^{\beta_{n-2}}NE^{\beta_{n-1}}E,\eea where $N$ and $E$ are steps
north and east, and exponentiation indicates repetition of a step.
We claim that $P$ and $Q$ are Dyck paths from $(0,0)$ to
$(n-4,n-4)$, and that $P$ never goes below $Q$. We define
$\bij(T)=(P,Q)$.

\begin{figure}[hbt]
\epsfig{file=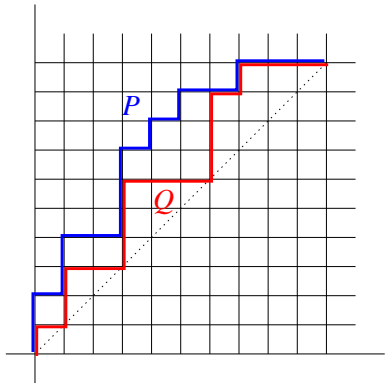,width=2in} \caption{\label{fig:pairdyck}
The pair $\bij(T)=(P,Q)$, where $T$ is the $2$-triangulation from
Figure~\ref{fig:big2triang}.}
\end{figure}

For example, if $T$ is the $2$-triangulation from
Figure~\ref{fig:big2triang}, we get from
Figure~\ref{fig:big2triang1}(i) that \bea\nn P=NNENNEENNNENENEENEEE,
\\ Q=NENNEENNNEEENNNENEEE. \nn\eea
These paths are drawn in Figure~\ref{fig:pairdyck}.

\ms

We claim that at each step of the coloring algorithm there is always
a cross to be colored red and a cross to be colored blue in the
appropriate blocks, so all crosses get colored at the end. We have
also stated that $P$ and $Q$ are non-crossing Dyck paths. Finally,
we claim that $\bij$ is in fact a bijection between $\TT_n$ and
$\DD_{n-4}$. We will justify these assertions in the next three
sections, by giving more insight on the bijection. The idea is to
construct isomorphic generating trees for the set of
$2$-triangulations and the set of pairs of non-crossing Dyck paths.
The natural isomorphism between the two generating trees determines
$\bij$.

\section{A generating tree for $2$-triangulations} \label{sec:2triang}

In this section we describe a generating tree where nodes at level
$\ell$ correspond to $2$-triangulations of an $(\ell+5)$-gon. The
root of the tree is the only $2$-triangulation of a pentagon, which
has no diagonals.

In the rest of this paper, when we refer to a $2$-triangulation we
will not consider the trivial diagonals. In particular, all
$2$-triangulations of an $n$-gon have $2n-10$ diagonals. The {\em
degree} of a vertex is the number of (nontrivial) diagonals that
have it as an endpoint. The degree of $a$ is denoted $\deg(a)$.

\subsection{The parent of a $2$-triangulation}
To describe the generating tree, we specify the parent of any given
$2$-triangulation of an $n$-gon, where $n\ge6$. For this purpose we
need a few simple lemmas.

\begin{lemma}\label{lemma:inside}
Let $T\in\TT_n$ be a $2$-triangulation containing the diagonal
$(a,b)$, with $a<b-3$. Then $T$ contains the diagonal $(a,b-1)$ or a
diagonal of the form $(a',b)$ with $a<a'\le b-3$.
\end{lemma}

\begin{proof}
Assume that $(a,b-1)$ is not in $T$. Then, since $T$ is a maximal
set of diagonals with no $3$-crossings, adding the diagonal
$(a,b-1)$ would create a $3$-crossing together with two diagonals in
$T$. But these two diagonals together with $(a,b)$ do not form a
$3$-crossing. This means that at least one of these two diagonals
crosses $(a,b-1)$ but not $(a,b)$. This can only happen if such a
diagonal is of the form $(a',b)$ with $a<a'\le b-3$.
\end{proof}

\begin{lemma}\label{lemma:short}
Let $T\in\TT_n$ be a $2$-triangulation containing the diagonal
$(a,b)$, with $a\le b-3$. Then there exists a vertex
$i\in\{a,\ldots,b-3\}$ such that $T$ contains the diagonal
$(i,i+3)$.
\end{lemma}

\begin{proof}
If follows easily by iterating Lemma~\ref{lemma:inside}.
\end{proof}

\begin{lemma}\label{lemma:noshort}
Assume that $n\ge6$, and consider the labels of the vertices to be
taken modulo~$n$ (for example, vertex $n+1$ would be vertex $1$).
Let $T\in\TT_n$ be a $2$-triangulation that does not contain the
diagonal $(a,a+3)$. Then the degrees of the vertices $a+1$ and $a+2$
are both nonzero.
\end{lemma}

\begin{proof}
Since $T$ is a maximal set of diagonals without no $3$-crossing,
adding the diagonal $(a,a+3)$ would create a $3$-crossing. This can
only happen if in $T$ there is a diagonal with endpoint $a+1$ and
another diagonal with endpoint $a+2$ that cross.
\end{proof}

\begin{lemma}\label{lemma:deg0}
Assume that $n\ge6$, and consider the labels of the vertices to be
taken modulo~$n$. Let $T\in\TT_n$ be a $2$-triangulation and let $a$
be a vertex whose degree is $0$. Then $T$ contains the diagonals
$(a-2,a+1)$ and $(a-1,a+2)$.
\end{lemma}

\begin{proof}
If $(a-2,a+1)$ was not in $T$, then by Lemma~\ref{lemma:noshort} the
degree of $a$ would be nonzero. Similarly if $(a-1,a+2)$ was not in
$T$.
\end{proof}

Now we can define the parent of any given $2$-triangulation. Let
$n\ge6$, and let $T$ be a $2$-triangulation of an $n$-gon. Let $r$
be the largest number with $1\le r\le n-3$ such that $T$ contains
the diagonal $(r,r+3)$. This number $r=r(T)$ will be called the {\em
corner} of $T$. Diagonals of the form $(i,i+3)$ will be called {\em
short diagonals}.

Let us note look at some useful properties of $T$. First, note that
$T$ does not contain any diagonals of the form $(a,b)$ with $r<a\le
b-3\le n-3$, since otherwise, by Lemma~\ref{lemma:short}, there
would be a short diagonal contradicting the choice of $r$. In
particular, $T$ has no diagonals of the form $(r+1,b)$ or $(r+2,b)$
with $r+4\le b\le n$. We also have that $r\ge2$. Indeed, if $r=1$
then all the diagonals would have to be of the form $(1,b)$, but
there can only be $n-5$ such diagonals, which is half of the number
needed in a $2$-triangulation. There are three possibilities for the
degrees of the vertices $r+1$ and $r+2$.

If the degree of $r+2$ is zero, then by Lemma~\ref{lemma:deg0} the
diagonal $(r+1,r+4)$ belongs to $T$. In this case we have
necessarily that $n=r+3$, in order not to contradict the choice of
$r$, and this diagonal is in fact $(1,r+1)$.

If the degree of $r+1$ is zero, again by Lemma~\ref{lemma:deg0} we
have that $(r-1,r+2)$ belongs to $T$.

If the degrees of $r+1$ and $r+2$ are both nonzero, let $i$ be the
smallest index so that the diagonal $(i,r+1)$ belongs to $T$, and
let $j$ be the largest index so that the diagonal $(j,r+2)$ belongs
to $T$. By the previous reasoning, we know that $i,j<r$. It is also
clear that $j\le i$, since otherwise the diagonals $(i,r+1)$,
$(j,r+2)$ and $(r,r+3)$ would form a $3$-crossing. We claim that in
fact $i=j$. Indeed, by Lemma~\ref{lemma:inside} applied to the
diagonal $(j,r+2)$, we have that either $(j,r+1)$ belongs to $T$, in
which case $i\le j$ by the choice of $i$, or there is a diagonal in
$T$ of the form $(j',r+2)$ with $j<j'$, which would contradict the
choice of $j$.

With these properties in mind, we define the parent of $T$ in the
generating tree to be the $2$-triangulation $p(T)\in\TT_{n-1}$
obtained as follows:
\begin{itemize}
\item Delete the diagonal $(r,r+3)$ from $T$ (recall that $r:=\max\{a:\ 1\le a\le n-3,\ (r,r+3)\in T\}$).
\item If $\deg(r+1)=0$, delete the diagonal
$(r-1,r+2)$;\\
if $\deg(r+2)=0$ (in which case $r=n-3$), delete the
diagonal $(1,r+1)$;\\
if $\deg(r+1)>0$ and $\deg(r+2)>0$, delete the diagonal $(j,r+2)$,
where $j:=\max\{a:\ 1\le a<r,\ (a,r+2)\in T\}$ (in this case we also
have $j=\min\{a:\ 1\le a<r,\ (a,r+1)\in T\}$).
\item Contract the side $(r+1,r+2)$ of the polygon (that is, move
all the diagonals from $r+2$ to $r+1$, delete the vertex $r+2$, and
decrease by one the labels of the vertices $b>r+2$).
\end{itemize}
It is clear that $p(T)$ contains no $3$-crossings, because it has
been obtained from $T$ by deleting diagonals. Also, by the above
reasoning, $p(T)$ has exactly $2$ diagonals less than $T$.
Therefore, $p(T)$ is a $2$-triangulation of an $(n-1)$-gon.

\ms

\begin{figure}[hbt]
\epsfig{file=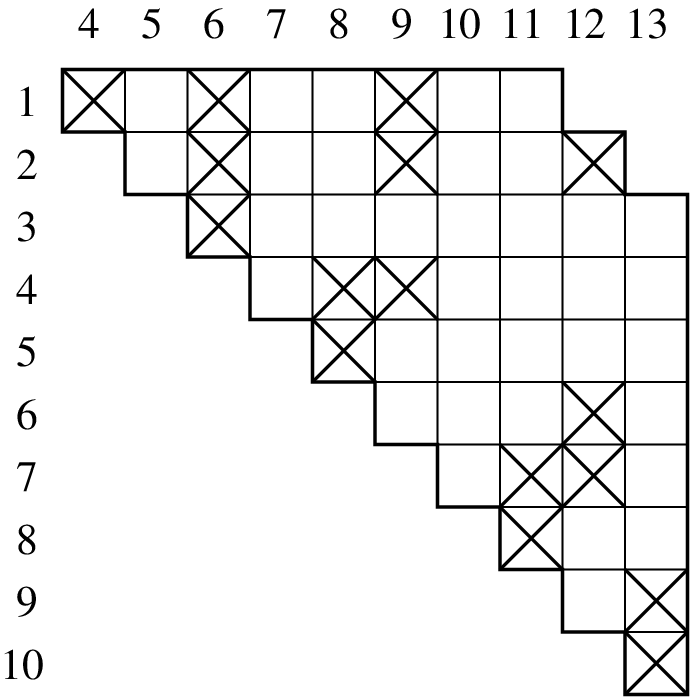,width=1.8in}\qquad\epsfig{file=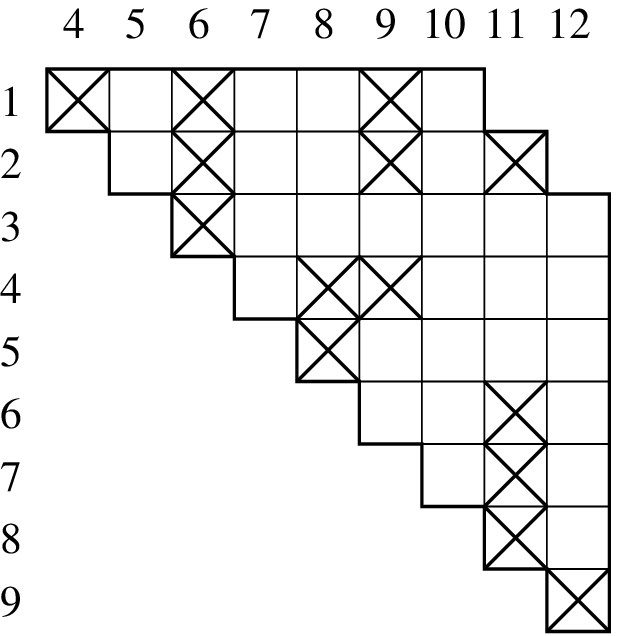,width=1.65in}
\qquad\epsfig{file=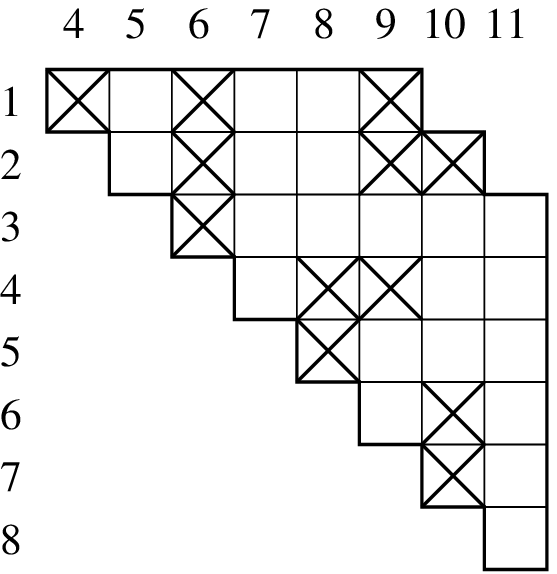,width=1.5in}
\caption{\label{fig:parent2tr} From left to right, the parent, the
grandparent, and the great grandparent of the $2$-triangulation from
Figure~\ref{fig:big2triang}.}
\end{figure}

It will be convenient to give an equivalent description of $p(T)$ in
terms of diagrams of $2$-triangulations. Consider the representation
of $T$ as a subset of $\La_n$. Next we describe how the diagram of
$p(T)$ as a subset of $\La_{n-1}$ is obtained from it. Observe that
if $r$ is the corner of $T$, then the diagram of $T$ has no crosses
below row $r$, because crosses in squares $(a,b)$ with $r<a\le
b-3\le n-3$ would contradict the choice of $r$, by
Lemma~\ref{lemma:short}.
 To obtain the diagram of $p(T)$, first delete all the squares
$(a,a+3)$ for $a=r-1,r,\ldots,n-3$. (Note that aside from $(r,r+3)$,
the only square among these where there may be a cross is
$(r-1,r+2)$, and if this cross is present, then column $r+1$ is
empty.) Next we merge columns $r+1$ and $r+2$. We do this so that
the new merged column, which will be the new column $r+1$, has a
cross in those rows where either the old column $r+1$ or $r+2$ (or
both) had a cross. (Note that there is at most one row where both
columns had a cross.) This yields the diagram of $p(T)$ as a subset
of $\La_{n-1}$. For example, if $T$ is the $2$-triangulation from
Figure~\ref{fig:big2triang}, then $p(T)$, $p(p(T))$ and $p(p(p(T)))$
are shown in Figure~\ref{fig:parent2tr}.

Note that in the bijection $\bij$ defined in Section~\ref{sec:bij},
the iterated step that merges blocks $r-2$ and $r-1$ consists
precisely in moving up one level in this generating tree of
$2$-triangulations. At each iteration, if $n'-3$ is the current
number of blocks, this indicates that we have moved up in the tree
to a $2$-triangulation $T'$ of a $n'$-gon. Then, for $1\le a\le b\le
n'-3$, a cross in row $a$ and block $b$ indicates that the diagonal
$(a,b+3)$ is present in $T'$. The largest $r$ such that there is a
cross in row $r$ and block $r$ is the precisely the corner of $T'$.
Merging blocks $r-2$ and $r-1$ in the original diagram is equivalent
to merging columns $r+1$ and $r+2$ in $T'$.

\subsection{The children of a $2$-triangulation}
Even though the generating tree is already completely specified by
the above subsection, it will be useful to characterize the children
of a given $2$-triangu\-lation $T\in\TT_n$ in the tree. By
definition, the children are all those elements $\wh{T}\in\TT_{n+1}$
such that $p(\wh{T})=T$. Again, let $r\in\{1,2,\ldots,n-3\}$ be the
corner of $T$. Equivalently, $r$ is the largest index of a nonempty
row in the diagram of $T$. Note that for any child $\wh{T}$ of $T$,
if $\hat{r}$ is the corner of $\wh{T}$, one must have $\hat{r}\ge
r$. It is not hard to check that all the children of $T$ are
obtained in the following way:
\begin{itemize}
\item Choose a number $u\in\{r,\ldots,n-2\}$.
\item Add one to the labels of the columns $j$ with $u+2\le j\le n$.
\item Add the square $(u,u+3)$ with a cross in it, and add empty
squares $(j,j+3)$ for $j=u+1,\ldots,n-2$.
\item Split column $u+1$ into two columns labeled $u+1$ and $u+2$ as
follows: \begin{enumerate}
\item Let $(a_1,u+1),\ldots,(a_h,u+1)$ be the crosses in column $u+1$ (assume that $a_1>\cdots>a_h$).
Choose a number $i\in\{0,1,\ldots,h\}$. If $u=n-2$, there is an
additional available choice $i=h+1$; if this is chosen, skip to (5)
below.
\item Leave the crosses $(a_1,u+1),\ldots,(a_i,u+1)$ in column
$u+1$.
\item Add a cross in position $(a_i,u+2)$ if $i>0$, or in position
$(u-1,u+2)$ if $i=0$.
\item Move the crosses $(a_{i+1},u+1),\ldots,(a_h,u+1)$ to
$(a_{i+1},u+2),\ldots,(a_h,u+2)$.
\item In the special case that $u=n-2$ and that $i=h+1$ has been chosen, column $u+1$ is split by leaving all the
crosses $(a_1,u+1),\ldots,(a_h,u+1)$ in it, adding a new cross
$(1,u+1)$, and leaving column $u+2$ empty.
\end{enumerate}
\end{itemize}

\begin{figure}[hbt]
\epsfig{file=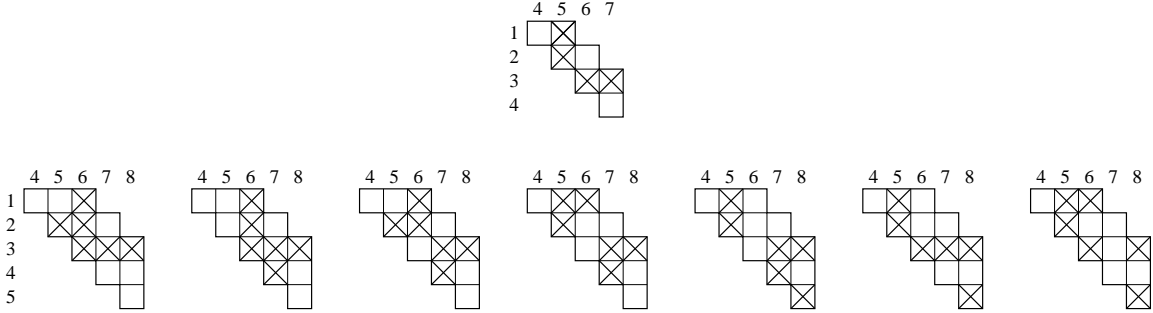,width=6in}
\caption{\label{fig:children2tr} A $2$-triangulation of an heptagon
and its 7 children in the generating tree.}
\end{figure}

Each choice of $u$ and $i$ gives rise to a different child of $T$.
Note that each choice of $u$ generates those children with
$\hat{r}=u$. Figure~\ref{fig:children2tr} shows a $2$-triangulation
and its seven children, of which one is obtained with $u=3$, three
with $u=4$, and three with $u=5$. It follows from the above
characterization that the total number of children of $T$ is
$$(h_{r+1}+1)+(h_{r+2}+1)+\cdots+(h_{n-1}+1)+1=h_{r+1}+h_{r+2}+\cdots+h_{n-1}+n-r,$$
where, for $r<j<n$, $h_j$ is the number of crosses in column $j$ of
the diagram of $T$. This observation allows us to easily describe
the generating tree for $2$-triangulations by labeling the nodes
with the list of numbers $(h_{r+1},\ldots,h_{n-1})$. For each chosen
$u\in\{r,\ldots,n-2\}$, the $h_{u+1}$ crosses in column $u+1$ can be
split into two columns for each choice of $i$. We have proved the
following result.

\begin{prop}\label{prop:gt2triang} The generating tree described above for the set $\TT$
is isomorphic to the tree with root labeled $(0,0)$ and with
generating rule
$$(d_1,d_2,\ldots,d_s)\longrightarrow \begin{array}{l}
\{(i,d_j-i+1,d_{j+1}+1,d_{j+2},\ldots,d_s)\ :\ 1\le j\le s-1,\ 0\le
i\le d_j \}\\ \cup\ \{(i,d_s-i+1)\ :\ 0\le i\le d_s+1\}.
\end{array}$$
\end{prop}

\begin{figure}[hbt]
\epsfig{file=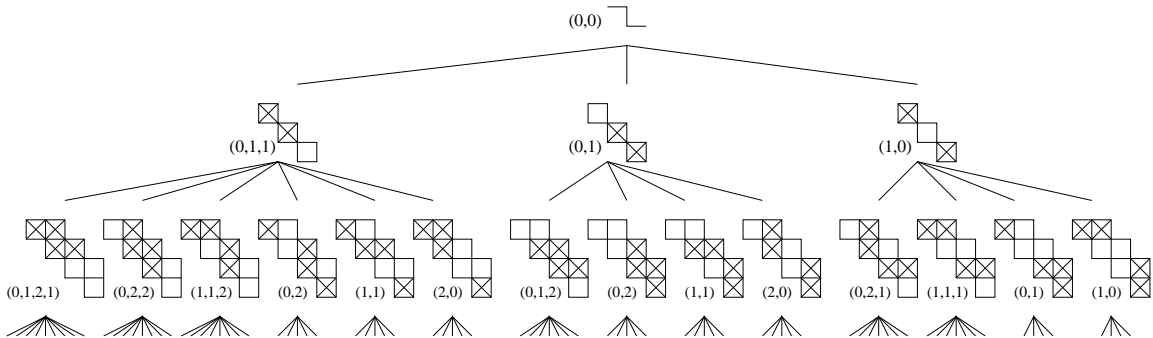,width=6in}
\caption{\label{fig:gentree2triang} The first levels of the
generating tree for $2$-triangulations.}
\end{figure}

For example, the children of a node labeled $(0,1,3,2)$ have labels
$(0,1,2,3,2)$, $(0,2,4,2)$, $(1,1,4,2)$, $(0,4,3)$, $(1,3,3)$,
$(2,2,3)$, $(3,1,3)$, $(0,3)$, $(1,2)$, $(2,1)$, and $(3,0)$. In
Figure~\ref{fig:children2tr}, the parent has label $(0,2,1)$ and the
children, from left to right, are labeled $(0,1,3,1)$, $(0,3,2)$,
$(1,2,2)$, $(2,1,2)$, $(0,2)$, $(1,1)$, and $(2,0)$. The first
levels of the generating tree for $\TT$ with their labels are drawn
in Figure~\ref{fig:gentree2triang}.

\section{A generating tree for pairs of non-crossing Dyck
paths} \label{sec:Dyck}

In this section we define a generating tree for $\DD$, where nodes
at level $\ell$ correspond to pairs of Dyck paths of size $\ell+1$
such that the first never goes below the second, and we show that it
is isomorphic to the generating tree from
Proposition~\ref{prop:gt2triang}. The root of our tree is the pair
$(P,Q)$, where $P=Q=NE$.

Every Dyck path $P$ of size $m$ can be expressed uniquely as
$$P=NE^{p_m}NE^{p_{m-1}}\cdots NE^{p_2}NE^{p_1}E$$
for some nonnegative integers $p_i$. The sequence
$(p_1,p_2,\ldots,p_m)$ determines the path, and it must satisfy
$p_1+p_2+\cdots+p_t\ge t-1$ for all $1\le t\le m$, and
$p_1+p_2+\cdots+p_m=m-1$. Given a pair $(P,Q)\in\DD_m$, we will
write $P$ as above, and $Q$ as
$$Q=NE^{q_m}NE^{q_{m-1}}\cdots NE^{q_2}NE^{q_1}E.$$
We set $p_{m+2}=p_{m+1}=q_{m+1}=0$ by convention. It will be
convenient to encode the pair $(P,Q)$ by the matrix
$$[P,Q]:=\left[\begin{array}{cccccccc}
p_{m+2}&p_{m+1}&p_m&p_{m-1}&\cdots&p_3&p_2&p_1 \\
q_{m+1}&q_m&q_{m-1}&q_{m-2}&\cdots&q_2&q_1&0\end{array}\right].$$
The leftmost column has zero entries, so it is superfluous, but it
will make the notation easier later on. The condition that $P$ never
goes below $Q$ is equivalent to the fact that for any
$t\in\{1,\ldots,m\}$, $p_1+p_2+\cdots+p_t\ge q_1+q_2+\cdots+q_t$. We
will write $p_j(P,Q)$ and $q_j(P,Q)$ when we want to emphasize that
these are parameters of the pair $(P,Q)$. We define
$$s=s(P,Q)=\min\{j\ge2 : p_jq_j=0\}.$$
Note that $2\le s\le m+1$. For example, the encoding of the pair
$(P,Q)$ of paths in Figure~\ref{fig:pairdyck} is
$$[P,Q]=\left[\begin{array}{cccccccccccc}
0&0&0&1&0&2&0&0&1&1&2&2 \\
0&1&0&2&0&0&3&0&0&1&2&0\end{array}\right]$$ and $s(P,Q)=3$.

The parent of $(P,Q)$ in the generating tree is defined to be the
pair $(P',Q')\in\DD_{m-1}$ whose encoding is
$$[P',Q']:=\left[\begin{array}{ccccccccccc}
p_{m+2}&p_{m+1}&p_m&\cdots&p_{s+2}&p_{s+1}+p_{s}&p_{s-1}-1&p_{s-2}&\cdots&p_2&p_1 \\
q_{m+1}&q_m&q_{m-1}&\cdots&q_{s+1}&q_{s}+q_{s-1}-1&q_{s-2}&q_{s-3}&\cdots&q_1&0\end{array}\right].$$
Note that in the case that $s=m+1$, both $[P,Q]$ and $[P',Q']$ have
the form \beq\label{eq:degen_paths}\left[\begin{array}{cccccccc}
0&0&1&1&\cdots&1&1&0 \\
0&1&1&1&\cdots&1&0&0\end{array}\right].\eeq

If we let $s'=s(P',Q')$, then it is clear from the definitions that
$s'\ge s-1$. Finally, observe that $P'$ never goes below $Q'$ since,
by the choice of $s$, we must have $p_{s}=0$ or $q_{s}=0$. For
example, the parent of the pair of Dyck paths drawn in
Figure~\ref{fig:pairdyck} is
$$[P',Q']=\left[\begin{array}{cccccccccccc}
0&0&0&1&0&2&0&0&2&1&2 \\
0&1&0&2&0&0&3&0&0&2&0\end{array}\right].$$

\ms

The above description completely specifies the generating tree for
$\DD$. As in the case of $2$-triangulations, it will be useful to
characterize the children of the pair $(P,Q)\in\DD_m$. Let
$p_j,q_j$, for $j=1,\ldots,m$, and $s$ be defined as above. The
children are the pairs $(\wh{P},\wh{Q})\in\DD_{m+1}$ whose parent
$((\wh{P})',(\wh{Q})')$ obtained using the above construction is
again $(P,Q)$. Note that if $\hat{s}=s(\wh{P},\wh{Q})$, then
$\hat{s}\le s+1$. It is easy to check that the children of $(P,Q)$
are the pairs $(\wh{P},\wh{Q})$ obtained in the following way.

\begin{itemize}
\item Choose a number $t\in\{1,2,\ldots,s\}$.
\item The following are the encodings of the children of
$(P,Q)$: \beq\label{eq:ch1}
[\wh{P},\wh{Q}]=\left[\begin{array}{ccccccccccc}
p_{m+2}&p_{m+1}&\cdots&p_{t+2}&p_{t+1}-i&i&p_{t}+1&p_{t-1}&\cdots&p_2&p_1 \\
q_{m+1}&q_{m}&\cdots&q_{t+1}&0&q_{t}+1&q_{t-1}&q_{t-2}&\cdots&q_1&0\end{array}\right]\eeq
for each $i\in\{1,\ldots,p_{t+1}\}$,
\beq\label{eq:ch2}[\wh{P},\wh{Q}]=\left[\begin{array}{ccccccccccc}
p_{m+2}&p_{m+1}&\cdots&p_{t+2}&p_{t+1}&0&p_{t}+1&p_{t-1}&\cdots&p_2&p_1 \\
q_{m+1}&q_{m}&\cdots&q_{t+1}&0&q_{t}+1&q_{t-1}&q_{t-2}&\cdots&q_1&0\end{array}\right],\eeq
and
\beq\label{eq:ch3}[\wh{P},\wh{Q}]=\left[\begin{array}{ccccccccccc}
p_{m+2}&p_{m+1}&\cdots&p_{t+2}&p_{t+1}&0&p_{t}+1&p_{t-1}&\cdots&p_2&p_1 \\
q_{m+1}&q_{m}&\cdots&q_{t+1}&j&q_{t}-j+1&q_{t-1}&q_{t-2}&\cdots&q_1&0\end{array}\right]\eeq
for each $j\in\{1,\ldots,q_{t}\}$ if $t\ge2$, or
$j\in\{1,\ldots,q_{t}+1\}$ if $t=1$.
\end{itemize}

\begin{figure}[hbt]
\epsfig{file=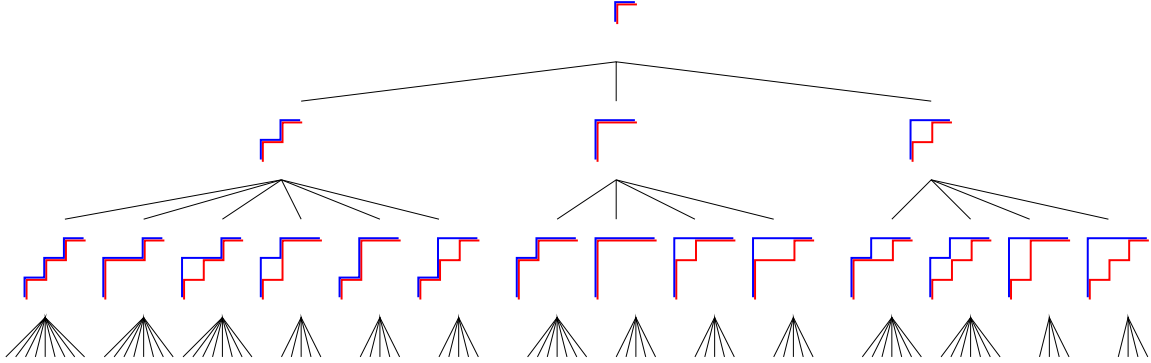,width=6in}
\caption{\label{fig:gentreepairsdyck} The first levels of the
generating tree for $2$-triangulations.}
\end{figure}

Note that each choice of $t$ generates the children with
$\hat{s}=t+1$. This is why when the column of $[P,Q]$ with entries
$p_{t+1}$ and $q_{t}$ is split into two columns, say
$\mathrm{col}_\mathrm{left}$ and $\mathrm{col}_\mathrm{right}$,
either the upper entry of $\mathrm{col}_\mathrm{right}$ or the lower
entry of $\mathrm{col}_\mathrm{left}$ has to be $0$. The first
levels of the generating tree for $\DD$ are drawn in
Figure~\ref{fig:gentreepairsdyck}.

\section{Why is $\bij$ a bijection?}\label{sec:why}

In this section we proof that $\bij$ is indeed a bijection. We start
by showing that the generating tree for pairs of non-crossing Dyck
paths from the previous section is the same as the one we
constructed for $2$-triangulations.

\begin{theorem}\label{th:isogt}
The generating tree for $\TT$ given in Section~\ref{sec:2triang} is
isomorphic to the generating tree for $\DD$ given in
Section~\ref{sec:Dyck}.
\end{theorem}

\begin{proof}
For our generating tree for $2$-triangulations,
Proposition~\ref{prop:gt2triang} gives a simple description of the
generating rule, with an appropriate labeling of the nodes. All we
need to show is that we can assign labels to pairs of non-crossing
Dyck paths so that our tree for $\DD$ obeys the same generating
rule.

Given a pair $(P,Q)\in\DD_m$, let
$p_1,p_2,\ldots,p_{m+1},p_{m+2},q_1,q_2,\ldots,q_m,q_{m+1}$, and
$s=s(P,Q)$ be defined as in Section~\ref{sec:Dyck}. We define the
label associated to the corresponding node of the tree to be
$$(p_{s+1}+q_{s},p_{s}+q_{s-1},\ldots,p_2+q_1).$$
Note that the root is labeled $(0,0)$.

For each node $(P,Q)$ in the tree for $\DD$, each choice of
$t\in\{1,2,\ldots,s\}$ yields children $(\wh{P},\wh{Q})$ with
$\hat{s}=s(\wh{P},\wh{Q})=t+1$. If $t\ge2$, then the number of
children generated by a particular choice of $t$ is
$p_{t+1}+q_{t}+1$, and their labels, according to (\ref{eq:ch1}),
(\ref{eq:ch2}), (\ref{eq:ch3}), and the above definition, are
$$\left\{\begin{array}{lcll}
(p_{t+1}-i,&\hspace{-2mm}q_{t}+i+1,&\hspace{-2mm}p_{t}+q_{t-1}+1,p_{t-1}+q_{t-2},\ldots,p_2+q_1)
& \mbox{for each $i\in\{1,\ldots,p_{t+1}\}$,} \\
(\ \
p_{t+1},&\hspace{-2mm}q_{t}+1,&\hspace{-2mm}p_{t}+q_{t-1}+1,p_{t-1}+q_{t-2},\ldots,p_2+q_1),
& \mbox{and} \\
(p_{t+1}+j,&\hspace{-2mm}q_{t}-j+1,&\hspace{-2mm}p_{t}+q_{t-1}+1,p_{t-1}+q_{t-2},\ldots,p_2+q_1)
& \mbox{for each $j\in\{1,\ldots,q_{t}\}$,}
\end{array}\right.$$ or equivalently,
$$(l,\ p_{t+1}+q_{t}-l+1,\ p_{t}+q_{t-1}+1,p_{t-1}+q_{t-2},\ldots,p_2+q_1)
\quad \mbox{for each $l\in\{1,\ldots,p_{t+1}+q_{t}\}$.}$$ Similarly,
the choice $t=1$ generates $p_{2}+q_{1}+2$ children, whose labels
are
$$(l,\ p_{2}+q_{1}-l+1)
\quad \mbox{for each $l\in\{1,\ldots,p_{2}+q_{1}+1\}$.}$$ This is
clearly equivalent to the generating rule from
Proposition~\ref{prop:gt2triang}, so the theorem is proved.
\end{proof}

Note that in the generating trees in the above proof, the labels of
the children of any particular node are all different. This uniquely
determines an isomorphism of the generating trees, which in turn
naturally induces a bijection $\wt\Psi$ between $2$-triangulations
of an $n$-gon and pairs of Dyck paths of size $n-4$ so that the
first never goes below the second. Let us analyze some properties of
this bijection. Consider a $2$-triangulation $T\in\TT_n$ and its
corresponding pair $\wt\Psi(T)=(P,Q)\in\DD_{n-4}$. Then, the
parameter $r$ in $T$ and the parameter $s$ in $(P,Q)$ are related by
$r+s=n-1$. The value of $u\in\{s,\ldots,n-2\}$ chosen to generate a
child of $T$ and the value of $t\in\{1,\ldots,s\}$ chosen to
generate a child of $(P,Q)$ are related by $u+t=n-1$. Also, if
$h_j$, for $j=r+1,\ldots,n-1$, is defined to be the number of
crosses in column $j$ of the diagram of $T$, and $p_j,q_j$, for
$j=1,\ldots,n-2$, are defined as above, then the label
$(d_1,\ldots,d_s)$ of the nodes corresponding to $T$ and $(P,Q)$ is
\beq\label{eq:label}
(d_1,\ldots,d_s)=(h_{r+1},\ldots,h_{n-1})=(p_{s+1}+q_{s},p_{s}+q_{s-1},\ldots,p_2+q_1).\eeq

Given a $2$-triangulation $T\in\TT_n$, in order to compute
$\wt\Psi(T)$ we find the path in the tree from the node
corresponding to $T$ to the root, keeping track of the labels of the
nodes encountered along the path. Then, starting from the root
$(NE,NE)$ in the generating tree for $\DD$, these labels determine
how to descend in the tree level by level, until we end with a pair
$(P,Q)$ of Dyck paths of size $n-4$, which is $\wt\Psi(T)$ by
definition. In a similar way we can compute the inverse
$\wt\Psi^{-1}((P,Q))$, where $(P,Q)\in\DD_m$.

For example, consider $T\in\TT_{14}$ to be the $2$-triangulation
represented in Figure~\ref{fig:big2triang}. Its corner is $r=10$,
and the label of the corresponding node in the tree for
$2$-triangulations is $(1,2,4)$, since those are the numbers of
crosses in columns $11$, $12$ and $13$, respectively. Its parent,
shown in the left of Figure~\ref{fig:parent2tr}, has $r=10$ and
label $(2,3)$. Its grandparent, drawn in the middle of
Figure~\ref{fig:parent2tr}, has $r=9$ and label $(0,4)$. Its great
grandparent has $r=7$ and label $(2,3,3)$. If we continue going up
in the generating tree, the next labels that we get are $(0,4,2)$,
$(0,3,3,1)$, $(0,1,2,2,1)$, $(0,1,2,1)$, $(0,1,1)$, and $(0,0)$, the
last one being the label of the root. To obtain $\wt\Psi(T)$, we
start with the root of the tree for $\DD$, whose encoding is
$\left[\begin{array}{ccc} 0&0&0
\\ 0&0&0\end{array}\right]$. Of its three children, the one with label
$(0,1,1)$ is generated by rule (\ref{eq:ch2}) with $t=2$, and its
encoding is $\left[\begin{array}{cccc} 0&0&1&0
\\ 0&1&0&0\end{array}\right]$. The next node down the tree with label
$(0,1,2,1)$ is encoded by $\left[\begin{array}{ccccc} 0&0&1&1&0
\\ 0&1&1&0&0\end{array}\right]$. Its child with label $(0,1,2,2,1)$ is
$\left[\begin{array}{cccccc} 0&0&1&1&1&0
\\ 0&1&1&1&0&0\end{array}\right]$. Rule (\ref{eq:ch1}) with $t=3$ and
$i=1$ generates the next node $\left[\begin{array}{ccccccc}
0&0&0&1&2&1&0
\\ 0&1&0&2&1&0&0\end{array}\right]$, with label $(0,3,3,1)$. Again, rule
(\ref{eq:ch1}) with $t=2$ and $i=2$ generates its child $\left[\begin{array}{cccccccc} 0&0&0&1&0&2&2&0 \\
0&1&0&2&0&2&0&0\end{array}\right]$, with label $(0,4,2)$. Rule
(\ref{eq:ch2}) with $t=2$ generates the next node
$\left[\begin{array}{ccccccccc} 0&0&0&1&0&2&0&3&0
\\ 0&1&0&2&0&0&3&0&0\end{array}\right]$, with label $(2,3,3)$. Its child with label
$(0,4)$ is generated using rule (\ref{eq:ch1}) with $t=1$, and it is
$\left[\begin{array}{cccccccccc} 0&0&0&1&0&2&0&0&3&1 \\
0&1&0&2&0&0&3&0&1&0\end{array}\right]$. Following the path down
according to the labels we got, we obtain pairs of Dyck paths whose
encodings are
$\left[\begin{array}{ccccccccccc} 0&0&0&1&0&2&0&0&2&1&2 \\
0&1&0&2&0&0&3&0&0&2&0\end{array}\right]$,
 and $\left[\begin{array}{cccccccccccc} 0&0&0&1&0&2&0&0&1&1&2&2 \\
0&1&0&2&0&0&3&0&0&1&2&0\end{array}\right]$. The last one is by
definition the encoding of $\wt\Psi(T)$, which is the pair in
Figure~\ref{fig:pairdyck}.

We claim that $\wt\Psi$ is precisely the bijection $\bij$ defined in
Section~\ref{sec:bij}. The description that we gave of $\bij$ is
nonrecursive, although implicitly it also computes the path to the
root in the generating tree for $\TT$. To justify this claim we use
the following lemma.

\begin{lemma}\label{lemma:columns}
Fix $n\ge5$. Let $T\in\TT_n$, and let
$(P,Q)=\wt\Psi(T)\in\DD_{n-4}$. For $4\le j\le n$, let $h_j$ be the
number of crosses in column $j$ of the representation of $T$ as a
subset of $\Lambda_n$. For $1\le j\le n-4$, let $p_j=p_j(P,Q)$ and
$q_j=q_j(P,Q)$.
Then,
$$(h_4,h_5,\ldots,h_{n-1},h_{n})=(q_{n-4},p_{n-4}+q_{n-5},\ldots,p_{2}+q_{1},p_1).$$
\end{lemma}

\begin{proof}
First notice that equation~(\ref{eq:label}) shows that the lemma
holds for the rightmost $s$ components not including the last one,
where $s=s(P,Q)$.

We prove the lemma by induction on $n$. For $n=5$, the empty
$2$-triangulation has $h_4=h_5=0$, and the pair of Dyck paths of
size one has $p_1=q_1=0$. Assume now that $n\ge6$ and the result
holds for $n-1$. Given $T\in\TT_n$, let $T'=p(T)\in\TT_{n-1}$ be its
parent, and let $(P',Q')=\wt\Psi(T')$. For $4\le j\le n-1$, let
$h'_j$ be the number of crosses in column $j$ of the representation
of $T'$ as a subset of $\Lambda_{n-1}$. For $1\le j\le n-5$, let
$p'_j=p_j(P',Q')$ and $q'_j=q_j(P',Q')$, and let $p'_{n-4}=q'_0=0$.
By the induction hypothesis,
$(h'_4,h'_5,\ldots,h'_{n-2},h'_{n-1})=(q'_{n-5},p'_{n-5}+q'_{n-6},\ldots,p'_{2}+q'_{1},p'_1)$.

Let $r$ be the corner of $T$, as usual. Let us first assume that
$r\ge3$. It follows from the rule that describes the children of
$T'$ that \beq\label{eq:h}
(h_4,h_5,\ldots,h_{n-1},h_{n})=(h'_4,\ldots,h'_{r},i,h'_{r+1}-i+1,h'_{r+2}+1,h'_{r+3},\ldots,h'_{n-1})\eeq
for some $0\le i\le h'_{r+1}$ if $r\le n-4$, or $0\le i\le
h'_{r+1}+1$ if $r=n-3$. Similarly, using that $s=n-1-r$, rules
(\ref{eq:ch1}), (\ref{eq:ch2}) and (\ref{eq:ch3}) describing the
children of $(P',Q')$ imply that
\bea\label{eq:ab} (q_{n-4},p_{n-4}+q_{n-5},\ldots,p_{2}+q_{1},p_1)=\hspace{6.5cm}\\
\nn
(p'_{n-4}+q'_{n-5},\ldots,p'_{s+1}+q'_{s},i,p'_{s}+q'_{s-1}-i+1,p'_{s-1}+q'_{s-2}+1,
p'_{s-2}+q'_{s-3},\ldots,p'_1+q'_0).\hspace{-12mm}\eea  We claim
that the value of $i$ has to be the same in (\ref{eq:h}) and
(\ref{eq:ab}). This is because by the definition of $\wt\Psi$, the
label of $T$ has to agree with the label of $(P,Q)$; but these
labels are given by the the rightmost $s$ components, not including
the last one, of (\ref{eq:h}) and (\ref{eq:ab}) respectively, and
the first entry is $i$ in both labels. It follows that (\ref{eq:h})
and (\ref{eq:ab}) coincide, so the lemma holds.

In the special case $r=2$, all the crosses in the diagram of $T$
have to be in the first two rows, and we have that
$(h_4,h_5,\ldots,h_{n-1},h_{n})=(1,2,2,\ldots,2,1,0)$. In this case,
$s=n-1-r=n-3$, and $[P,Q]$ has the form given in
(\ref{eq:degen_paths}), so we have that
$(q_{n-4},p_{n-4}+q_{n-5},\ldots,p_{2}+q_{1},p_1)=
(1,2,2,\ldots,2,1,0)$ as well.
\end{proof}

As an example of the fact stated in this lemma, take $T$ to be the
$2$-triangulation from Figure~\ref{fig:big2triang}, for which we
have seen that $\wt\Psi(T)$ is then the pair $(P,Q)$ of Dyck paths
drawn in Figure~\ref{fig:pairdyck}. In this case we have that
$$(h_3,h_4,\ldots,h_{13})=(1,0,3,0,2,3,0,1,2,4,2).$$
On the other hand, $$[P,Q]=\left[\begin{array}{cccccccccccc}
0&0&0&1&0&2&0&0&1&1&2&2 \\
0&1&0&2&0&0&3&0&0&1&2&0\end{array}\right],$$ so
$(q_{10},p_{10}+q_{9},\ldots,p_{2}+q_{1},p_1)=(0,1,0,3,0,2,3,0,1,2,4,2)$
as well.

\begin{figure}[hbt]
\epsfig{file=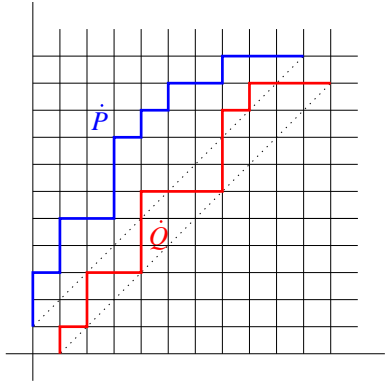,width=2in}
\caption{\label{fig:pairdyckmoved} The paths $\dot{P}$ and
$\dot{Q}$, where $(P,Q)$ are drawn in Figure~\ref{fig:pairdyck}.}
\end{figure}

A convenient way to represent a pair $(P,Q)\in\DD_m$ is to shift the
paths slightly, drawing $P$ as a path from $(0,1)$ to $(m,m+1)$,
which we call $\dot{P}$, and $Q$ as a path from $(1,0)$ to
$(m+1,m)$, which we call $\dot{Q}$ (see
Figure~\ref{fig:pairdyckmoved}). The fact that $P$ does not go below
$Q$ is equivalent to the fact that $\dot{P}$ and $\dot{Q}$ do not
intersect. In the drawing of $\dot{P}$ and $\dot{Q}$, the number of
east steps with ordinate $j$ is then $p_{m-j+2}+q_{m-j+1}$ for
$j=1,\ldots,m-1$; $p_2+q_1+1$ for $j=m$; and $p_1+1$ for $j=m+1$.

Lemma~\ref{lemma:columns} states that if $T\in\TT_n$ and
$(P,Q)=\wt\Psi(T)$, then the number $h_j$ of crosses in column $j$
of $\La_n$ equals the number of east steps with ordinate $j-3$ in
the drawing of $(\dot{P},\dot{Q})$ (except when $j$ equals $n-1$ or
$n$, where these numbers are off by $1$). This explains why in the
definition of $\bij$ we considered the number of crosses in each
column of $\La_n$ to determine where to put the east steps in $P$
and $Q$. It remains to see how many of these $h_j$ east steps belong
to $P$ and how many belong to $Q$, that is, how to split $h_j$ into
$p_{n-j+1}+q_{n-j}$.

In the definition of $\bij$, this is given by coloring the crosses
red and blue. To determine how to color the crosses, let us analyze
now the encodings of the children of a fixed pair $(P,Q)\in\DD_m$.
Let $s=s(P,Q)$, and let $t\in\{1,\ldots,s\}$ be the parameter chosen
to generate a particular child of $(P,Q)$.

Rules (\ref{eq:ch1}), (\ref{eq:ch2}) and (\ref{eq:ch3}) show that
the $(t+1)$-st column from the right of $[P,Q]$ (the one with
entries $p_{t+1}$ and $q_{t}$) is split into two columns, and then a
$1$ is added to the bottom entry of the new right column and to the
top entry of the column immediately to the right of it. Thus, the
first of these $1$'s contributes to $\wh{Q}$, and the second one to
$\wh{P}$. This explains why in the iterated step of the description
of $\bij$, a cross in block $r$ is colored blue (contributing to the
upper path) and a cross in the block to the left of it is colored
red (contributing to the lower path).

Now, the blocks encountered in this iterated step are, in general,
sets of adjacent columns of $\La_n$ that have been merged when going
up in the tree for $2$-triangulations. So, how do we know, among all
the crosses in a block, which is the one that has to be colored red
(or blue)? The key observation is that whenever a column $\bt{c} $p_{t+1}$\\
$q_{t}$\et$ of $[P,Q]$ is split into two columns, according to rules
(\ref{eq:ch1}), (\ref{eq:ch2}), and (\ref{eq:ch3}), the upper entry
gets split only if the lower entry moves entirely to the right
column, and the lower entry gets split only if the upper entry moves
entirely to the left column. This means that in a block that
consists of merged columns, a cross that contributes to the lower
(resp. upper) path will always come from the rightmost (resp.
leftmost) possible column among the merged ones. So, when a cross in
a block that consists of merged columns needs to be colored red
(resp. blue), we must always color the rightmost (resp. leftmost)
uncolored cross in the block.

Note that in case of a tie, that is, if there is more than one
rightmost (or leftmost) uncolored cross, it does not matter which
one we color. This is because the construction of $(P,Q)=\bij(T)$
only takes into account the {\em number} of red and blue crosses in
each column of the diagram, but not which particular crosses have
each color.

\section{Generalization to $k$-triangulations}\label{sec:k}

The natural question at this point is whether one can give a similar
bijection between $k$-triangu\-lations of an $n$-gon and $k$-tuples
$(P_1,P_2,\ldots,P_k)$ of Dyck paths of size $n-2k$ such that each
$P_i$ never goes below $P_{i+1}$, for $k\ge3$. While we have not
succeeded in finding such a bijection, some of the ideas in our
construction for $k=2$ generalize to arbitrary $k$. In this section
we show that it is possible to construct an analogous generating
tree for $k$-triangulations.

\subsection{A generating tree for $k$-triangulations}

Fix an integer $k\ge2$. Next we describe a generating tree where
nodes at level $\ell$ correspond to $k$-triangulations of an
$(\ell+2k+1)$-gon. We ignore trivial diagonals, so all
$k$-triangulations of an $n$-gon have $k(n-2k-1)$ diagonals. The
root of the tree is the empty $k$-triangulation of a $(2k+1)$-gon.

The lemmas in Section~\ref{sec:2triang} have an immediate
generalization to arbitrary $k$. We will only use two of them.

\begin{lemma}\label{lemma:insidek}
Let $T\in\Tk_n$ be a $k$-triangulation containing the diagonal
$(a,b)$, with $a<b-k-1$. Then $T$ contains the diagonal $(a,b-1)$ or
a diagonal of the form $(a',b)$ with $a<a'\le b-k-1$.
\end{lemma}

\begin{proof}
Assume that $(a,b-1)$ is not in $T$. Then, since $T$ is a maximal
set of diagonals with no $(k+1)$-crossings, adding the diagonal
$(a,b-1)$ would create a $(k+1)$-crossing together with $k$
diagonals in $T$. But these $k$ diagonals together with $(a,b)$ do
not form a $(k+1)$-crossing. This means that at least one of these
$k$ diagonals crosses $(a,b-1)$ but not $(a,b)$. This can only
happen if such a diagonal is of the form $(a',b)$ with $a<a'\le
b-k-1$.
\end{proof}

\begin{lemma}\label{lemma:shortk}
Let $T\in\Tk_n$ be a $k$-triangulation containing the diagonal
$(a,b)$, with $a\le b-k-1$. Then there exists a vertex
$i\in\{a,\ldots,b-k-1\}$ such that $T$ contains the diagonal
$(i,i+k+1)$.
\end{lemma}

Lemma~\ref{lemma:shortk} follows easily by iteration of
Lemma~\ref{lemma:insidek}.





Diagonals of the form $(a,a+k+1)$ are called {\em short diagonals}.
Let $n\ge2k+2$, and let $T$ be a $k$-triangulation of an $n$-gon. To
define the parent of $T$ we will need some definitions. Let $r$ be
the largest number with $1\le r\le n-k-1$ such that $T$ contains the
short diagonal $(r,r+k+1)$. We call $r$ the {\em corner} of $T$.
Note that $T$ does not contain any diagonals of the form $(a,b)$
with $r<a\le b-k-1\le n-k-1$, since otherwise, by
Lemma~\ref{lemma:shortk}, there would be a short diagonal
contradicting the choice of $r$. So, the diagram of $T$ has no
crosses below row $r$. Note that in particular we have $r\ge k$,
since each $a\le r$ can be an endpoint of at most $n-2k-1$
diagonals, compared to the $k(n-2k-1)$ needed in a
$k$-triangulation.

\begin{figure}[hbt]
\epsfig{file=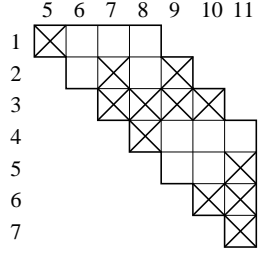,width=1.3in} \caption{\label{fig:3triang}
The diagram of a $3$-triangulation of an $11$-gon.}
\end{figure}

For $i=1,2,\ldots,k-1$, let $$A_i:=\{a:(a,r+i)\in T\} \cup
\{r+i-k\}.$$ Let $a_1:=\min A_1$, and for $i=2,\ldots,k-1$, let
$$a_i:=\min\{a\in A_i : a>a_{i-1}\}.$$ For example, in the
$3$-triangulation from Figure~\ref{fig:2triang}, $r=7$, $a_1=3$, and
$a_2=6$. The following property of $T$ will be crucial to define its
parent.

\begin{lemma}\label{lemma:a_i} Let $i\in\{1,2,\ldots,k\}$, and let $a_i$ be defined as above.
Then, either $(a_i,r+i+1)\in T$ or $(a_i,r+i+1)$ is a trivial
diagonal.
\end{lemma}

\begin{proof}
First notice that if $a\in A_i$, then $a\le r+i-k$. This is because
the diagram of $T$ has no crosses below row $r$, so all diagonals
incident to $r+i$ are represented by crosses in column $r+i$, whose
lowest square is in row $r+i-k-1$.

We start with the case $i=1$. If the square $(a_1,r+2)$ falls
outside of $\La^{(k)}_n$, then $(a_1,r+2)$ is a trivial diagonal and
we are done. Otherwise, let us assume for contradiction that
$(a_1,r+2)\not\in T$. Since $T$ is a maximal set of diagonals with
no $(k+1)$-crossings, this means that if we added $(a_1,r+2)$ to
$T$, it would form a $(k+1)$-crossing together with $k$ diagonals in
$T$, none of which corresponds in the diagram to a cross below row
$r$ (since there are no such crosses). By the definition of $a_1$,
none of these diagonals can correspond to a cross in column $r+1$.
Therefore, if in this $(k+1)$-crossing we replace $(a_1,r+2)$ with
$(a_1,r+1)$, we obtain a $(k+1)$-crossing containing $(a_1,r+1)$,
which contradicts the fact that $T$ is a $k$-triangulation.

For $i>1$ the reasoning is very similar. In this case, we assume for
contradiction that $(a_i,r+i+1)\not\in T$ and that it is not a
trivial diagonal. Then, maximality of the set $T$ implies that
adding $(a_i,r+i+1)$ would create a $(k+1)$-crossing $C$, together
with $k$ diagonals in $T$. By the definition of
$a_1,a_2,\ldots,a_i$, there must be at least one among the columns
$r+1,r+2,\ldots,r+i$ which has no diagonals belonging to $C$. Let
$r+j$ be the rightmost such column. Then, if for each
$l=j,j+1,\ldots,i$ we replace the element in $C$ in column $r+l+1$
with $(a_l,r+l)$, we still obtain a $(k+1)$-crossing. But the fact
that the diagonal $(a_i,r+i)$ is part of a $(k+1)$-crossing is a
contradiction, since either $(a_i,r+i)\in T$ or $(a_i,r+i)$ is a
trivial diagonal.
\end{proof}

An additional property of $T$ is that column $r+k$ of its diagram
has no crosses below row $a_{k-1}$. This is because if there was
such a cross, then it would form a $(k+1)$-crossing together with
diagonals $(a_1,r+1),\ldots,(a_{k-1},r+k-1)$, and $(r,r+k+1)$, all
of which belong to $T$ or are trivial diagonals.

\ms

Consider now the representation of $T$ as a subset of $\La^{(k)}_n$.
We define the parent of $T$ in the generating tree to be the
$k$-triangulation $p(T)\in\Tk_{n-1}$ whose diagram, as a subset of
$\La^{(k)}_{n-1}$, is obtained from the diagram of $T$ as follows.

\bit
\item Delete the squares $(a,a+k+1)$ for $a=r,r+2,\ldots,n-k-1$. (Note
that only the first one of such squares contains a cross.)
\item For each $i=1,2,\ldots k-1$: \bit
\item Keep all the crosses of the form $(a,r+i)$ with $a\ge a_i$ in
column $r+i$.
\item Move all the crosses of the form $(a,r+i+1)$ with $a<a_i$ from
column $r+i+1$ to column $r+i$, and delete the cross $(a_i,r+i+1)$
if it is in $T$. \eit
\item Delete column $r+k$ (which at this point is empty, by the observation following
Lemma~\ref{lemma:a_i}), and move all the columns to the right of it
one position to the left. If $r>n-2k$, delete also the squares
$(a,n-k-1+a)$ for $a=1,2,\ldots,r+2k-n$. \eit

\begin{figure}[hbt]
\epsfig{file=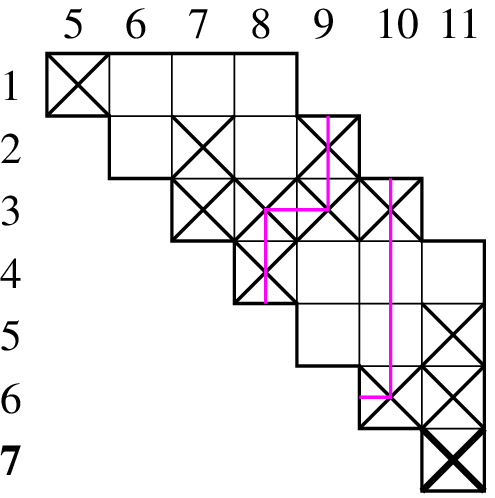,width=1.3in}\hspace{1cm}\bt{c}$\longrightarrow$\\
\\ \\ \\ \\
\et\hspace{1cm}\epsfig{file=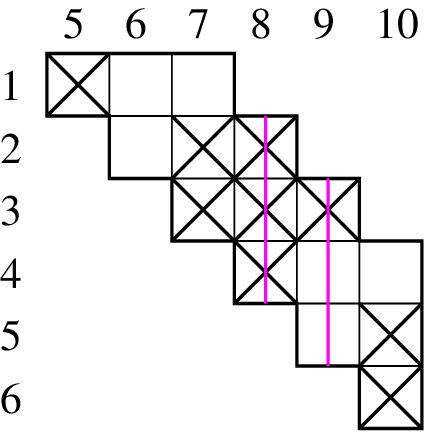,width=1.1in}\vspace{-7mm}

\epsfig{file=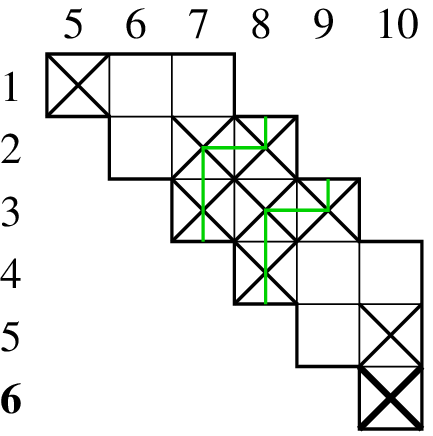,width=1.1in}\hspace{1cm}\bt{c}$\longrightarrow$\\
\\ \\ \\ \et\hspace{1cm}\epsfig{file=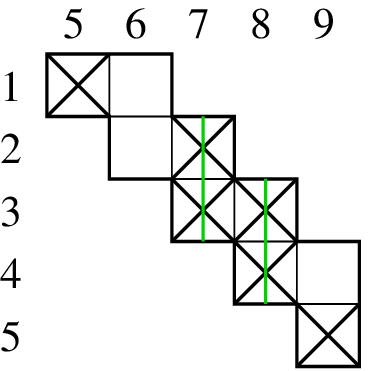,width=.93in}\vspace{-8mm}
\caption{\label{fig:parent3tr} The parent (top) and the grandparent
(bottom) of the $3$-triangulation from Figure~\ref{fig:3triang}.}
\end{figure}

This yields the diagram of $p(T)$ as a subset of $\La^{(k)}_{n-1}$.
For example, if $T$ is the $3$-triangulation from
Figure~\ref{fig:3triang}, then $p(T)$ and $p(p(T))$ are shown in
Figure~\ref{fig:parent3tr}. Note that in $p(T)$, $r=6$, $a_1=2$, and
$a_2=3$.

\ms

We next characterize the children of a given $k$-triangulation
$T\in\Tk_n$ in the generating tree. By definition, the children are
all those elements $\wh{T}\in\Tk_{n+1}$ such that $p(\wh{T})=T$.
Again, let $r\in\{1,2,\ldots,n-k-1\}$ be the corner of $T$. Note
that for any child $\wh{T}$, if $\hat{r}$ is the corner of $\wh{T}$,
then $\hat{r}\ge r$. All the children of $T$ are obtained in the
following way:
\begin{itemize}
\item Choose a number $u\in\{r,\ldots,n-k\}$.
\item Add one to the labels of the columns $j$ with $u+k\le j\le n$, and add an empty column labeled $u+k$.
\item Add the square $(u,u+k+1)$ with a cross in it, and add empty
squares $(j,j+k+1)$ for $j=u+1,\ldots,n-k$. If $u>n-2k$, add also
empty squares $(j,n-k+j)$ for $j=1,\ldots,u+2k-n$.
\item For $i=1,\ldots,k-1$, let $B_i:=\{b:(b,u+i)\in T\} \cup
\{u+i-k\}$. If $u=n-k$, add also the element $i$ to $B_i$, for each
$i$.
\item For each $i=1,\ldots,k-1$, choose a number $b_i\in B_i$,
so that $b_1<b_2<\cdots<b_{k-1}$.
\item For each $i=k-1,k-2,\ldots,1$, add a cross $(b_i,u+i+1)$ (except if $b_i=i$, in which case
we add the cross $(b_i,u+i)$ instead), and move all the crosses of
the form $(b,u+i)$ with $b<b_i$ from column $u+i$ to column $u+i+1$.
\end{itemize}

\begin{figure}[hbt]
\epsfig{file=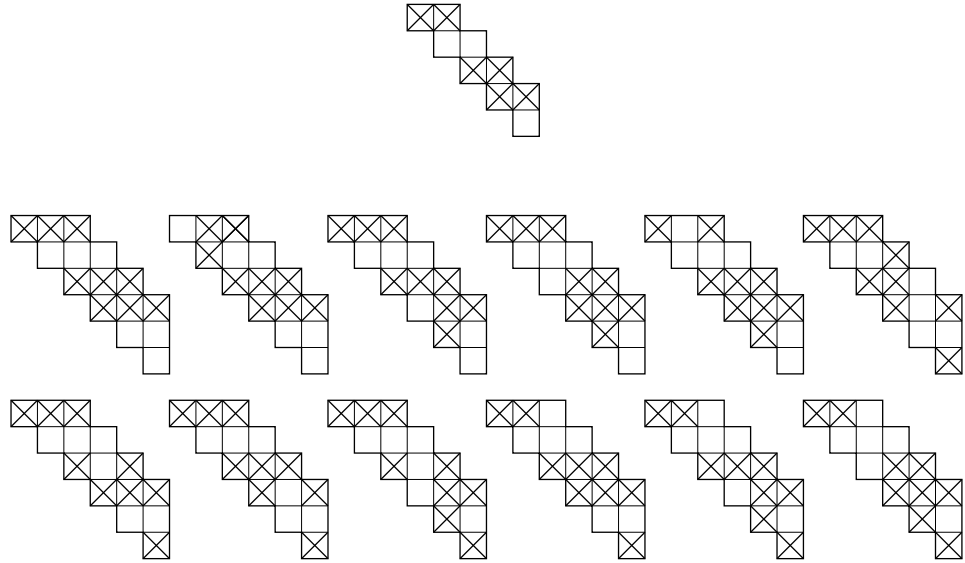,width=5in}
\caption{\label{fig:children3tr} A $3$-triangulation of a $9$-gon
and its 12 children in the generating tree.}
\end{figure}

Each choice of $u$ and $b_1,b_2,\ldots,b_{k-1}$ gives rise to a
different child of $T$. Note that each choice of $u$ generates those
children with $\hat{r}=u$. Figure~\ref{fig:children3tr} shows a
$3$-triangulation and its twelve children, of which two are obtained
with $u=4$, three with $u=5$, and seven with $u=6$. An important
difference between the case $k=2$ and the case $k\ge3$ is that, in
the latter, the number of children of a $k$-triangulation depends
not only on the number of crosses in the columns of its diagram but
also on the relative position of the crosses in different columns
(this is caused by the condition $b_1<b_2<\cdots<b_{k-1}$). As a
consequence, there is no obvious way to associate simple labels to
each node of the generating tree, as we did for $k=2$. This is an
obstacle when trying to construct a generating tree for $k$-tuples
of non-crossing Dyck paths isomorphic to the one that we have given
for $\Tk$.

\ms

\subsection*{Acknowledgements}
I am grateful to Marc Noy for interesting conversations, and for
suggesting the idea of trying a recursive approach find a bijection.
I also thank Richard Stanley and Peter Winkler for helpful
discussions.

\ms

\end{document}